\documentclass[11pt]{amsart}
\usepackage{geometry}                
\usepackage{graphicx}
\usepackage{amssymb}
\usepackage{epstopdf}
\usepackage[czech]{babel}
\title{Is There Ethics in Mathematics?}
\author{Roman Kossak}

\begin{document}
\maketitle
Intrigued by Borovik's  remarks on the Cambridge University Ethics in Mathematics Project
Mathematics Project (EiM) in his recent article {\it The Tool/Weapon Duality of Mathematics}, revisited [Selected Passages from Correspondence with Friends 13 no. 2 (2025), 15–31.  http://bit.ly/4l5LO5z], I decided to find out more about the work of EiM.  There is much to read.  It seems that a good place to start a discussion with the  document  titled {\it Four Levels of Ethical Engagement} [EiM Discussion Paper 1/2018 University of Cambridge Ethics in Mathematics Project, https://www.ethics-in-mathematics.com/assets/dp/18\_1.pdf] by Maurice Chiodo and Piers Bursill-Hall. {\bf Unless indicated otherwise, all quotes in my text are from this paper}. In the introduction, the authors state

\begin{quote} Nobody questions that doctors and lawyers need to have some ethical training, guidelines, and a sensitivity to difficult issues their practice may encounter. We assert the premise here that the same is true for mathematicians.[...]

This premise is far from universally accepted by mathematicians, in our experience. Many
mathematicians do not feel that there is ethics in mathematics, nor that there is any need for
subject-specific training in ethics for mathematicians.
\end{quote}
The authors  identify, classify, and give examples of fraudulent and nefarious applications of mathematics and point out that the ethical dilemmas  they give rise to are a call for a much more engaged
response from the mathematical community. However, one big obstacle they notice is that

\begin{quote} Many [mathematicians] feel that mathematics is value-free and that
mathematics is inherently extra-ethical, or that mathematics itself does not have ethical
implications or consequences. [...]  

In our view this is somewhat akin to gun manufacturers and arms merchants arguing “Guns don’t
kill people; people kill people,” and therefore that they have no ethical, moral or legal
responsibility in the killing of people using their guns. [...]

We have noticed that most mathematicians are implicitly {\it taught} to be, and think at, this level.
They deny that ethical issues are raised by doing mathematics and argue that any ethical
consequences that follow from their mathematics do not lie within the proper domain of the
working mathematician. They always point to someone else as having responsibility.
\end{quote} 

Clearly, something needs to be done, but first we need to find more about what makes us and our colleagues so insensitive to ethical  issues.  Is it because mathematics attracts people with such an attitude?
Is it because mathematics makes people  insensitive? Is it because of  lack of proper instruction and guidance in education of future mathematicians?  A discussion of this takes most of the 25 pages of the article, but first  the authors introduce five levels of awareness and engagement in the
ethics in mathematics  to classify specific attitudes and behaviors of mathematicians. The levels are not hierarchical, but there is a lowest, the zero-level, which is described as   ``[...] to think that there is no ethics in mathematics, or that ethical consequences are not the responsibility of the mathematician".  The other levels---for the descriptions of which  I use the titles of the corresponding sections of the article---are: (1) Realising there are ethical issues inherent in mathematics; (2) Doing something: speaking out to other mathematicians; (3) Taking a seat at the tables of power; and (4) Calling out the bad mathematics of others. 

There is much to discuss here, but I will concentrate just on Level 0, and this is because it is the level I found myself at.  I do care about ethical conduct in one's personal life and in professional career  and I am very open to discussing ethical issues with anyone. I even sat at some tables of small power, but I cannot move to Level 1, because I do not believe that there are issues inherent in mathematics that would not be covered by much more general ethics considerations.  I will try to explain why, but let me first summarize  some conclusions I came to  while reading Borovik's and Chiodo and Bursill-Hall's articles. 
\begin{itemize} 
\item The picture of a “typical mathematician” that Chiodo and Bursill-Hall  paint is cartoonish and false.
\item  Ethics is a one of the hardest  parts of philosophy and one cannot seriously discuss any particular issues concerning ethics without at least acknowledging some of the general difficulties. 
\item One cannot learn how to act ethically  by taking a course in ethics. No amount of learned wisdom can prevent one from making ethically questionably decisions (case study: Heidegger).
\item However, one learns how to behave ethically by observing how other people behave (role models).
\end{itemize}
\section{Who are those mathematicians?}
The next two quotes seem to me essential to understand  Chiodo and Bursill-Hall's motivation for their project.
\begin{quote}
Whatever it is that mathematicians might do, we need to realise that any discussion of this
requires an analysis of how mathematicians behave. We believe this is an extremely important
condition, and one that mathematicians generally do not want to take into account.
Mathematicians are a self-selected and well-defined group of people: you know who you would
call a mathematician and who you would not. We need to understand better what the natural,
inherent characteristics of the members of the community, or of the profession, are. How do those
personality and intellectual characteristics function in various places in the social matrix of
academia or as part of the mathematical engine room of a corporation?  [...]

Many professional, working mathematicians are exactly the publicity-shy people who are not
going to enjoy the limelight, and do not enjoy controversy or public confrontation. They may be vicious, harsh, and ruthless when arguing about mathematics over email, without face-to-face
social contact, but neither the heated debate nor the social subtleties of the committee room or
the public podium are the most comfortable place for them. Moreover, many mathematicians are
not culturally well trained or well prepared to engage in public argument of an ethical, political,
indeterminate nature and the reader might note that nothing in our training as mathematicians
addresses this. Of course, this is a generalisation, but it is probably a reasonable generalisation of
the norms of mathematicians’ behaviour.
\end{quote}

Do you know such mathematicians? I am sure you do, and so do I. Is this a ``reasonable generalisation of
the norms of mathematicians’ behaviour"? I doubt. And one more quote from  the final section of the paper:

\begin{quote} However, the conclusion we feel is the most important for further study of ethics in mathematics
is that we need the skills and knowledge of social scientists – sociologists, anthropologists, and
psychologists -- (i) to understand better the nature of mathematicians as individuals from the
point of view of their particular psychological traits, and (ii) to understand the community values
and behaviour, the ‘tribal’ behaviour of mathematicians. 
\end{quote}

I cannot bring any social science insights into this discussion, but I can offer something else. It seems to me interesting that my vision of who mathematicians are  diverges so much from Chiodo and Bursill-Hall's. I do not doubt that their observations are honest, and the conclusions  they draw from them are sincere, but why are they so different from mine? If indeed there are  inherent aspects of our profession that are the reason that most---or even just many---of our colleagues behave in the ways described in the article, why is it  that  during my long career in mathematics I formed a quite different opinion about who we are and what we do? 

I was admitted to the M.S. program in  mathematics at the University of Warsaw  in Poland in 1972.  It is important to say that it was a large program that at that time admitted three hundred students per year. It was run by a big department  divided into the institutes of Mathematics, Computer Science, and Mechanics. Each institute was divided into several sections. In the first  year  we all took required courses: introduction to  mathematics, introduction to computer science, analysis, and algebra.  Then we had to choose between mathematics and computers science. I chose mathematics.   In the second year I took required courses in more analysis,  differential equations, topology, and mathematical
logic. The level of those courses was comparable to those given in the Ph.D. programs in the US. In the third year,  one had to choose a topic to specialize in.  I chose foundations of mathematics---which in Poland meant set theory, proof theory, computability theory, and model theory---and for the next three years I was taking courses in those disciplines, combined with some supplementary ones such as functional analysis and elementary number theory. In the fourth and fifth year, we were expected  to attend research seminars as well. The program also included  two-semester courses in  philosophy and physics. I am writing all this to stress that in my formative years in mathematics I met a large number of mathematicians from many diverse disciplines. 

At that time, all exams had both written and oral components.  Oral exams were given by the professors teaching the course, and some of those professors were legends of the Polish mathematics. I can say with confidence  that we got to know all of them relatively well. I do not remember a single person that would fit the stereotype described by Chiodo and Bursill.  They were a  diverse group, but at the same time quite united by their attitude towards mathematics as a subject and by the understanding of its role in the society.   

I found the  following statement  in a popular Polish mathematical magazine for high school students Delta.
\begin{quote}
\dots I believe that mathematics is the structure of the world. Not a description of this
structure, but the structure itself. A mathematician can undoubtedly create very
strange objects and it may seem to him that he has departed far from
reality. This is only an appearance. If it is good mathematics, it will turn out sooner or later
to be a fragment of reality\dots . If the world were different, it would be different
mathematics.  
\end{quote}
 It was taken  from an article by the Polish mathematician Andrzej Lasota  (my translation). The editors of Delta considered  it significant enough to be  published as a separate short note.  To me, it represents well  what made mathematics  attractive for me and for many of my fellow students.

\section{Engagement}

The mathematicians I got to know during my student years were  all quite passionate about mathematics,   but there was also something else that was very characteristic. It was the level of their engagement in politics and social issues.  Chiodo ans Bursill-Hall  wrote:

\begin{quote} The stereotypical mathematician precisely eschews all those socially difficult, under-determined
and emotionally complex sorts of loci of power and power exchange. Mathematicians are good at
mathematics, but often not so good at committees and positions of power where social
negotiation and nuance rule.
\end{quote}
Fine, but what about all those non-stereotypical ones---those that may not form a majority, or even a very large part of the community, but who nevertheless shape the attitudes of others and influence their decisions?  I will give examples.

The 1970s and 80s were difficult times in Poland. The system was crumbling, and the opposition to it was growing. The mathematics department at the University of Warsaw was a hub of dissident activities. Many students and faculty were directly involved in opposition to the government.  The activities were clandestine, but visible to us in different forms. It all  started  with  student  protests against the anti-Zionist campaign began in 1967 by a faction of the ruling party, and continued later against the Soviet invasion of Czechoslovakia in 1968. Some participants were banned from the university, some were sentenced to prison. The mathematics department was under constant surveillance. Some secret service informants were  embedded as fake students, but it was a common knowledge who they were. 

In 1980, the Solidarity movement became an official resistance force until it was crushed by the imposition of the martial law in December of 1981.  This eventually led to the collapse of the system, but only in 1989.  Many mathematicians were active agents  of those developments.  I will mention two of them.

In high school, I learned basics of the mathematical logic from a book written by Wiktor Marek and Janusz Onyszkiewicz. Later, as a student,  I took a model theory course given by Onyszkiewicz.  Onyszkiewicz was one of the founders of the Masovian Region of the Solidarity trade union. For his political activity  he was interned for a year after the martial law was imposed in 1981. In 1986 he became the spokesman for the Solidarity union. I was already in the US then and I remember seeing him on television news delivering short and precise statements. In 1989,  he was one of the negotiators on behalf of the union at the round table talks with the government at which  agreements were made for a transition to democracy. After 1989, he served as Minister of Defense twice, and was very respected by the military for his professionalism.  In 2004, was elected to the European Parliament (and this is a very abridged version of his cv).

I learned set theory from Andrzej Zarach, who,  while working on his Ph.D.,  organized informal seminars for freshmen. Despite of his accomplishments in mathematics and excellent teaching skills,  he could not get an academic job in Warsaw  because of his political activity. While working at the Wroc{\l}aw University of Science and Technology, he became a founder  of the local Solidarity unit. He was a delegate to the first National Congress of Delegates in 1981 and a member of the Presidium of the Solidarity National Commission. In 1987, he organized the first underground information services in Poland, distributing information in electronic form. At the same time, he wrote articles for underground magazines, becoming a member of the editorial board of the magazine {\it Fighting Solidarity}. \footnote{An anecdote: In Wroc{\l}aw, Zarach offered a graduate course in proof theory. The course was not easy and quite technical and soon all students dropped out, except for one: a well-known police informant. Zarach continued teaching him until the end of the semester.} 

In my student days, I attended talks and had personal interactions with many foreign visitors. 
I learned about the kind of model theory in which I specialize in now from George Wilmers who was visiting  for a semester from the University of Manchester.  Wilmers was one of the organizers of the Bertrand Russell Memorial Logic Conference, which was held in Uldum in Denmark in 1971. It was a  counter-conference organized in opposition to one held in Cambridge with NATO funding which generated enormous controversy at the time.\footnote {Onyszkiewicz was also one of the organizers of the counter-conference.}  Here is the table of contents of Part I of the proceedings of the conference: Critical and Philosophical. 

\begin{itemize}
\item Bertrand Russell - a tribute arranged by Graham Priest.
\item Uldum 1971: the Bertrand Russell Memorrial Logic  Conference and the controversy surrounding it, by
Alan Slomson.
\item The New Universal Church by Alexandre Grothendieck (translated by Alan Slomson).
\item Ctitique or Grothendieck's paper, by Anders Kock.
\item A reply to Anders Kock, by Jan Waszkiewicz.
\item In defence or science, by L. D. Gasman.
\item Some remarks on current mathematical practice, by John Bell.
\item Artists and technicians, by Ross Skelton.
\item The predicament of intellectuals, by Aldo Ursini.
\item NATO's counter-revolutionary role: a case study by George Wilmers.
\item NATO and the myths or the Cold War, by Mark Cowling.
\item Note on the philosophy of mathematics, by Peter Clark.
\item A bedside reader's guide to the conventionalist philosophy of mathematics, by Graham Priest.
\item Some aspects or Wittgenstein's philosophy of mathematics, by Benedikt Peppinghaus.
\end{itemize}

The list of mathematical contributions includes papers by Max Dickmann, Per Lindstr\"om, Peter Aczel, Jeff Paris,  and Maurice Boffa---mathematicians I got to know personally. In my early years in mathematics, they were my role models. I learned mathematics  from them, but I also paid much attention to their high ethical standards.  There were many more who either visited Poland at that time or I had a chance to meet at conferences and who impressed me very much.  As I think now about those days, more and more names come to mind of  mathematicians from all over Europe and the US I got to know. All of them open mined and actively involved in politics at various levels. I will mention just two more names. The full list would be very long.

In Warsaw, we had strong connections to the logicians in what was then Czechoslovakia.  They were visiting often, except for their prominent leader Petr Vop\v{e}nka, who for political reasons was not allowed to travel abroad for conferences. Because of that I never met him, but I studied his work and was much impressed by it. In 1990, he was elected Vice Rector of Charles University and was Minister of Education 1990-1992. 

I first met Jouko V\"a\"an\"anen at a conference  in Poland in 1979.   In 1970s V\"a\"an\"anen joined the Finnish Union of Conscientious Objectors and was elected chair of the union. The group succeeded in relaxing some rules of  the law  concerning the obligatory military service.  Along with a distinguished  career in mathematics, for which he received the Magnus Ehrnrooth Foundation Prize in mathematics in  2024,    V\"a\"an\"anen has built a strong logic group at the University of Helsinki and was  a very successful  high level university administrator.

A particularly interesting event was the International Congress of Mathematicians that was planned to take place in Warsaw in 1982, but was cancelled because of the imposition of martial law in Poland. A decision was made to move the congress to 1983, but many  martial law restrictions were still in effect then and there were calls for boycott. I was among  those Polish mathematicians who were grateful  to the foreigners who had decided to participate, and even more grateful to some---including Ronald Graham---who in their lectures expressed solidarity with the suppressed Solidarity movement. They gave us a good lecture in ethics.

When I moved to the United States, I found out that there were differences  between mathematical academic life here and in Europe, but still at the City University of New York where I worked and elsewhere I found  many colleagues who showed no resemblance to the stereotypical mathematician. Some of them devoted their professional careers in equal parts to mathematics and to political activism and they were very effective in doing both. 

\section{Ethics and education}

Chiodo and Bursill end their paper with an appeal
\begin{quote}
Only social scientists have the skills and training
to analyse and explain to us how we work socially (if we can overcome {\it their} fear of
mathematicians and mathematics, and our notoriously superior attitudes towards their
disciplines). Only when we understand ourselves {\it in the new mathematical world of the 21st
century} will we be able to prepare ourselves and our students to become sensitive to, and
confront, both the ethical issues and the social responsibilities that come with our unique,
ubiquitous, world changing powers.

We cannot leave these issues to professional ethicists and philosophers. For the sake of
mathematics, and for the freedoms we enjoy, we must not just do nothing. No one else can do it,
so we must.
\end{quote}
To meet the challenge they propose to develop an ethics training designed specifically for our discipline. 
I am very skeptical, and one of the reasons is hidden in the quote above. The question is: How are we going  to deal with {\it their} fear of mathematics? The authors write:
\begin{quote} Mathematicians are notoriously unwilling or unable to explain their mathematics to the 
non-mathematician, and as a result there is often a vast disjunction between what a mathematician
knows her maths can do, what it means, how reliable or meaningful or probabilistic its results are
\dots and what the non-mathematical or semi-mathematical public believe the mathematics means.
\end{quote}
I absolutely agree that communication is a problem. I have been there, I have done that. During my career I had a chance to speak publicly at meetings of various levels of importance at my university. I am not a bad speaker and I have a reputation of being a good teacher. I am willing and able to explain mathematics to non-mathematicians; however, there are limits to what one can do, and that is because the non-mathematicians one talks to often know very  little about mathematics, and this is not necessarily their fault.\footnote{ A curious fact: In my academic discussions,  I learned that it is actually sometimes easier to talk to our colleagues in the humanities than in the disciplines that explicitly use mathematics. }

We are facing a long-lasting crisis in mathematical education in the US and  elsewhere. I have retired five years ago, but I hear from my colleagues that the three years of the pandemic and the massive move to online teaching made the already bad situation worse.  There is a vast literature about it, but for those who would like a quick insight into the roots and the complexity of the crisis I suggest {\it A Critique of the Structure of U.S. Elementary School Mathematics} by Liping Ma. [Notices of AMS, vol.60, No 10, 2013.]

 Chiodo and Bursill recommend a specific remedy:
\begin{quote}
Mathematicians often need to learn the specific skills required to
work with politicians, corporate management, and other non-scientists. These include engaging
in policy discussions, establishing and rationalising the objectives of their mathematical work,
and learning ways to communicate to a non-mathematical audience potential limitations,
overreach, and possible drawbacks.
\end{quote}
I do not really know what those specific skills might be. Policy discussions are by their very nature hard. Positions are biased, arguments can be misleading. What happens if one side brings in a mathematician who supports an idea, and  the other side a mathematician who opposes it? Whose competence should one trust?  Gian-Carlo Rota, in his {\it Ten rules for survival of a mathematics department} left us this advice that I think is very relevant:

\begin{quote} {\it Never wash your dirty linen in public.}

I know that you frequently (and loudly, if I may add) disagree with your colleagues about the relative value of fields of mathematics and about the talents of practicing mathematicians. All of us hold some of our colleagues in low esteem, and sometimes we cannot help ourselves from sharing these opinions with our fellow mathematicians.

When talking to your colleagues in other departments, however, these opinions should never be brought up. It is a mistake for you to think that you might thereby gather support against mathematicians you do not like. What your colleagues in other departments will do instead, after listening to you, is use your statements as proof of the weakness of the whole mathematics department, to increase their own departments' standing at the expense of mathematics.
\end{quote}
I think that it may not be a bad idea to include some training in lecturing (for everyone) and public speaking (for those who are thinking of doing it) in our Ph.D. programs in mathematics, but I do not think that it is essential. I also think that courses in history and philosophy of mathematics could be required for those who have never taken them before. In any case though, we should direct more effort  toward improving general mathematics education. 

A quick search shows that in Chiodo and Bursill's paper the word ``education" is used only once, and ``general education" never. This is a glaring omission. And this omission points to what seems to be missing from the whole argument---that only the broad social context gives the word ``ethics” meaning. Ethical behavior and its various patterns are not shaped by a particular profession or field of study. The ideals of higher education cannot be taught by sociologists, cannot be manufactured on request. There  must be a better way. Certainly, it’s worth pursuing.

\end{document}